\documentclass[11pt]{article}
\newtheorem{theorem}{Theorem}
\newtheorem{lemma}[theorem]{Lemma}
\newtheorem{prop}[theorem]{Proposition}
\newtheorem{coro}[theorem]{Corollary}

\newtheorem{remark}[theorem]{Remark}

 \newcommand{\nl}{\newline}
 \newcommand{\dist}{{\rm dist}}
 \newcommand{\N}{{\bf N}}
 
\newcommand{\R}{{\bf R}}

 \newcommand{\codim}{{\rm codim}}

\newcommand{\ia}{({\rm i})}
\newcommand{\ib}{({\rm ii})}
\newcommand{\ic}{({\rm iii})}
\newcommand{\id}{({\rm iv})}

\newcommand{\tc}{{ {\bf (*)}}}

\title{Best constants for higher-order Rellich inequalities in $L^p(\Omega)$}

\author{G. Barbatis}

\begin{document}

\maketitle

\

\begin{abstract}
We obtain a series improvement to higher-order $L^p$-Rellich inequalities on a Riemannian
manifold $M$. The improvement is shown to be sharp as each new term of the series is added.

\

\noindent {{\bf AMS Subject Classification:} 35J35 (35P15, 26D10)}\nl {{\bf Keywords:} Rellich inequality, best
constants, distance function }
\end{abstract}

\section{Introduction}

In the article \cite{DH} Davies and Hinz 
proved higher-order $L^p$ Rellich inequalities of the form
\begin{equation}
\int_{\R^N}\frac{|\Delta^mu|^p}{|x|^{\gamma}}dx \geq A(2m,\gamma)\int_{\R^N}\frac{|u|^p}{|x|^{\gamma+2mp}}dx\; ,
\label{sd}
\end{equation}
and
\begin{equation}
\int_{\R^N}\frac{|\nabla\Delta^mu|^p}{|x|^{\gamma}}dx \geq A(2m+1,\gamma)\int_{\R^N}\frac{|u|^p}{|x|^{\gamma+(2m+1)p}}dx\; ,
\label{sd1}
\end{equation}
for all $u\in C^{\infty}_c(\R^N\setminus \{0\})$
with the sharp value for the constants $A(2m,\gamma)$ and $A(2m+1,\gamma)$. Their approach uses some integral
inequalities involving $\Delta |x|^{\sigma}$
together with iteration, and is set initially in a Riemannian manifold context.
One of the aims of the present paper is to improve inequalities
(\ref{sd}) and (\ref{sd1}) by adding a sharp non-negative term at the respective right-hand sides.
In fact, this comes as a special -- and most important -- case of a more general
theorem where instead of $|x|$ we have a distance function $d(x)=\dist(x,K)$. 
Under a simple geometric assumption we establish an improved Rellich inequality of the form
\begin{equation}
\int_{\Omega}\frac{|\Delta^{m/2}u|^p}{d^{\gamma}}dx \geq A(m,\gamma)\int_{\Omega}\frac{|u|^p}{d^{\gamma+2mp}}dx
 + B(m,\gamma)\sum_{i=1}^{\infty}\int_{\Omega}V_i |u|^p \, dx\; , 
\label{chess}
\end{equation}
for all $u\in C^{\infty}_c(\Omega\setminus K)$,
where at each step we have an optimal function $V_i(x)$ and a
sharp constant $B(m,\gamma)$; see Theorem \ref{thm:meli} for the precise statement.
Here and below we interpret $|\Delta^{m/2}u|$ as $|\nabla \Delta^{(m-1)/2}u|$ when $m$ is odd.

Improved versions of Hardy or Rellich inequlities have attracted considerable attention recently and
especially for Hardy inequalities there is
a substantial literature; see, e.g., \cite{AE,BV,BFT2,GGM,T} and references therein.
The corresponding literature for Rellich inequalities is more restricted; see \cite{GGM,B,BT,TZ,T}.

As was the case in \cite{DH}, our results are formulated
in a Riemannian manifold context, but we note that they are also new in the Euclidean case.
We consider a Riemannian manifold $M$ of dimension $N\geq 2$, a domain $\Omega\subset M$,
a closed, piecewise smooth surface $K$ of codimension $k$, $1\leq k\leq N$, and the distance function $d(x):=\dist(x,K)$
which we assume to be bounded in $\Omega$. We note that this last assumption is only needed for the improved inequality
and not for the plain inequality where only the first term appears in the right-hand side of (\ref{chess});
to our knoweldge this is also new except in the case $M=\R^N$, $K=\{0\}$, studied in \cite{DH}.

We define recursively
\begin{eqnarray}
&& X_1(t) = (1- \log t)^{-1},\quad t\in (0,1],\nonumber \\
&& X_i(t) = X_{1}(X_{i-1}(t)), ~~~~~~~i=2,3,\ldots\, ,\quad  t\in (0,1].
\label{va}
\end{eqnarray}
These are iterated logarithmic functions that vanish at an
increasingly slow rate at $t=0$ and satisfy $X_i(1)=1$.

Given $m\in\N$ and $\gamma\geq 0$ we also define
\begin{eqnarray}
&&A'(m,\gamma)=\prod_{i=0}^{[(m-1)/2]}\Big( \frac{k-\gamma -(m-2i)p}{p}\Big)^p\; ,\nonumber\\
&& A''(m,\gamma)=\prod_{j=1}^{[m/2]}\Big( \frac{pk-k+\gamma +(m-2j)p}{p}\Big)^p\; , \nonumber\\
&& A(m,\gamma)=A'(m,\gamma)A''(m,\gamma)\; , \label{dan}\\
&& B(m,\gamma)=\frac{p-1}{2p}A(m,\gamma)
\Bigg\{  \sum_{i=0}^{[(m-1)/2]}
\Big( \frac{k-\gamma -(m-2i)p}{p}\Big)^{-2} + \nonumber\\
&&\qquad\qquad\qquad \qquad  +\sum_{j=1}^{[m/2]}
\Big( \frac{pk-k+\gamma +(m-2j)p}{p}\Big)^{-2}\Bigg\}\; .
\nonumber
\end{eqnarray}
Concerning the above definitions, we adopt the convention that empty sums are equal to zero and empty products are
equal to one; this of course refers to the sum or product over $j$ when $m=1$. To state our first theorem we
introduce the following technical hypothesis:
\[
\left(
\begin{array}{ll}
\gamma\neq \frac{3pk-8p^2-2k+6p}{4p-2}\; ,  & \mbox{if $m$ is even} \\
\gamma +p\neq \frac{3pk-8p^2-2k+6p}{4p-2} \; , & \mbox{if $m$ is odd, $m\geq 3$}
\end{array}
\right)
\quad \mbox{ or }\quad
p> \frac{13+\sqrt{105}}{4}\quad\quad \tc \; .
\]
We then have
\begin{theorem}
{\bf (improved Rellich inequality)} Let $m\in N$ and assume that $d(\cdot)$ is bounded in $\Omega$. Let $\gamma\geq 0$
be such that $k-\gamma-mp>0$ and suppose that $\tc$ is satisfied.
Assume moreover that
\begin{equation}
d\Delta d-k+1\geq 0 \quad , \quad \mbox{in $\Omega\setminus K$}
\label{cc}
\end{equation}
in the distributional sense. Then there exists a $D_0\geq \sup_{x\in\Omega}d(x)$ such that for any $D \geq  D_0$ there holds
\begin{eqnarray}
&&\int_{\Omega}d^{-\gamma} |\Delta^{m/2}u|^pdx \geq
A(m,\gamma)\int_{\Omega}d^{-\gamma-mp}|u|^p dx +\label{1.2}\\[0.2cm]
&&\qquad\qquad \qquad\qquad +B(m,\gamma)\sum_{i=1}^{\infty}\int_{\Omega}d^{-\gamma -mp}
 X_1^2X_2^2\ldots X_i^2 |u|^p dx,
\nonumber
\end{eqnarray}
for all  $u \in C^{\infty}_c(\Omega\setminus K)$, where $X_j=X_j(d(x)/D)$.
\label{thmb}
\end{theorem}
We present some examples where the geometric condition (\ref{cc}) is satisfied:

{\bf Example 1.} Suppose that $M=\R^N$ and that $K$ is affine. Then (\ref{cc}) is satisfied as an equality (this includes the case
where $K$ consists of a single point). 

{\bf Example 2.} Suppose that $M$ is a Cartan-Hadamard manifold, that is, a simply connected geodesically complete non-compact
manifold with non-positive sectional curvature. If $K=\{x_0\}$ (some point in $M$) then (\ref{cc}) is satisfied;
see \cite{SY}.

{\bf Example 3.} Suppose $M= M_1\times M_2$ where $M_1$ is a Cartan-Hadamard manifold of dimension $k$.
If $K=\{x_0\}\times M_2$ for some $x_0\in M_1$, then, in an obvious notation,
$d(x,y)(\Delta d)(x,y)=d(x,y)(\Delta d_1)(x)\geq d_1(x)(\Delta d_1)(x)\geq k-1$, so (\ref{cc}) is satisfied
for any $\Omega\subset M$.

Concerning the important case $M=\R^N$, $K=\partial\Omega$, condition
(\ref{cc}) is satisfied if $\Omega$ is the complement of a convex domain. This however is excluded from our theorem due to the assumptions
$k-\gamma-mp>0$, $\gamma\geq 0$
\footnote{We note here that the condition $\gamma\geq 0$ of Theorem \ref{thmb} can be weakened to $pk-k+\gamma>0$.
We have not persued this since our main goal is to obtain inequalities whose left-hand sides do not involve any weight
(and this only requires $\gamma\geq 0$; see the proof of Theorem 1 below).
We note however that if negative $\gamma$ were allowed, then condition $(*)$ would need
to be strengthened; see the comments in the beginning of the proof of Theorem \ref{thmb}.};
on the other hand, these conditions are not needed for Theorem \ref{thm:meli} below.
It should be noted that Rellich inequalities involving $\dist(x,\partial\Omega)$ present surprising difficulties when $p\neq 2$. In particular,
it is not known whether the inequality $\int_{\Omega}|\Delta u|^pdx \geq \{(p-1)(2p-1)/p^2\}^p\int_{\Omega}|u|^pd^{-2p}dx$
is valid when $\Omega$ is bounded and convex with a smooth boundary; see also \cite[Chapter 2]{T} for results in this direction.

In our second theorem we prove the optimality of the constants and exponents of Theorem \ref{thmb}. This is quite
technical and shall be established only in the case where $M=\R^N$ and $K$ is affine
(or, indeed, has an affine part, since the argument is local); we believe that extra effort should yield the
result in the general case, but we have not pursued this. This would require in particular estimates on the
behavior of higher-order derivatives of $d(x)$ near $K$; see \cite[Theorem 3.2]{AS}.
\begin{theorem}
{\bf (optimality)} Let $\Omega\subset\R^N$ and let $K$ be an affine hypersurface of codimension
$k\in \{1,\ldots,N\}$ such that $K\cap\Omega\neq \emptyset$. Assume that for some
$\gamma \in\R$, $D\geq \sup_{\Omega} d(x)$, $r\geq 1$ and some $\theta\in\R$, $C>0$ there holds
\begin{eqnarray*}
&&\int_{\Omega}d^{-\gamma} |\Delta^{m/2}u|^pdx \geq |A(m,\gamma)|\int_{\Omega}d^{-\gamma-mp}|u|^p dx + \\
&& \!\!\! +|B(m,\gamma)|\sum_{i=1}^{r-1}\int_{\Omega}d^{-\gamma-mp} X_1^2X_2^2\ldots X_i^2 |u|^p dx 
 +C\int_{\Omega}d^{-\gamma -mp}X_1^2X_2^2\ldots X_{r-1}^2X_r^{\theta} |u|^p dx \, ,
\end{eqnarray*}
for all $u\in C^{\infty}_c(\Omega)$, where $X_j=X_j(d(x)/D)$. Then
\begin{eqnarray*}
\ia &&  \theta \geq 2\; , \\
\ib && \mbox{if $\theta=2$ then $C\leq |B(m,\gamma)|$.}
\end{eqnarray*}
\label{thm:meli}
\end{theorem}
We note that $X_j(d(x)/D_1)/X_j(d(x)/D_2)\to 1$ as $d(x)\to 0$, for any $D_1,D_2\geq \sup_{\Omega}d$ and
in this sense the precise value of $D_0$ in Theorem \ref{thmb} does not affect the optimality of the theorem.
We also note that by a standard argument, if $k-\gamma-mp>0$ then the validity of (\ref{1.2}) for all $u\in C^{\infty}_c(\Omega\setminus K)$ implies
its validity for all $u\in C^{\infty}_c(\Omega)$.

The proof of Theorem 1 is given in Section 2 and uses some of the ideas of \cite{B} and, in particular, induction on $m$.
However, it is more technical due to $p\neq 2$ and the extra parameter $k$; 
moreover, the proof in \cite{B} uses one-dimensional arguments and depends on the Euclidean structure.
The proof of Theorem 2 is given in Section 3 and uses an appropriately chosen minimising sequence.

\section{The Rellich inequality}

Throughout the paper we shall repeatedly use the differentiation rule
\begin{equation}
 \frac{d}{dt} X_i^{\beta}(t)=\frac{\beta}{t}X_1(t)X_2(t)\ldots X_{i-1}(t)X_i^{\beta +1}(t),
\qquad i=1,2,\ldots,\quad \beta\in\R,
\label{diffrule}
\end{equation}
which is easily proved by induction on $i\in\N$. Let us define the functions
\[
\eta(t)=\sum_{i=1}^{\infty}X_1X_2\ldots X_i \; , \quad  \zeta(t)=\sum_{i=1}^{\infty}X_1^2X_2^2\ldots X_i^2 \; , 
\]
\[
\theta(t)=\sum_{i=1}^{\infty}\sum_{j=1}^iX_1^3 \ldots X_j^3 X_{j+1}^2\ldots X_i^2 \; ,
\]
(see \cite{B} for a detailed discussion of the convergence of these series). It follows from (\ref{diffrule}) that
\begin{equation}
\eta'(t)=\frac{\eta^2(t)+\zeta(t)}{2t}\; , \quad \zeta'(t)=\frac{2\theta(t)}{t}\; , \quad t\in (0,1)\; .
\label{nik}
\end{equation}

{\em Note.}
In the sequel we shall use the symbol $X_i$ for $X_i(t)$, $t\in (0,1]$, for
$X_i(d(x)/D)$, $x\in \Omega$, $D\geq \sup_{\Omega}d$, and also for $X_i(t/D)$, $t\in (0,\sup_{\Omega}d)$,
$D\geq\sup_{\Omega}d$. It will always be made explicit which
meaning is intended. The same also holds for the functions $\eta$, $\zeta$ and $\theta$.

For the sake of simplicity we work with real-valued functions, noting that with minor modifications the
proofs also work in the complex case.
As in \cite{DH}, the proof of Theorem \ref{thmb} uses iteration and for this we shall need to consider first the case
$m=2$. The following proposition has been obtained in \cite{BT} for $\gamma=\mu=0$. We set
\[
Q=\frac{(k-\gamma-2p)(pk-k+\gamma)}{p^2},
\]
\begin{prop}
Let $\gamma,\mu \geq 0$ be given and assume that $k-\gamma-2p>0$. Suppose that
\begin{equation}
d\Delta d-k+1 \geq 0 \quad , \quad \mbox{ in }\Omega\setminus K 
\label{woody}
\end{equation}
in the distributional sense and assume also that 
$(4p-2)\gamma\neq 3pk-8p^2-2k+6p$ or $p>(13+\sqrt{105})/4$.
Then there exists $D_0\geq \sup_{\Omega}d(x)$ such that for all $D\geq D_0$ and
all $\mu \geq 0$ there holds
\begin{eqnarray}
&&\quad \int_{\Omega}d^{-\gamma}(1+\mu\zeta)|\Delta u|^pdx \geq Q^p\int_{\Omega}d^{-\gamma-2p}|u|^pdx +
\label{ele}\\
&& +\left( \frac{p-1}{2p}Q^p\bigg\{ \Big(\frac{k-\gamma-2p}{p}\Big)^{-2} 
+\Big(\frac{pk-k+\gamma}{p}\Big)^{-2}\Big\}
+ Q^p\mu\right)\!\!\int_{\Omega}d^{-\gamma-2p}\zeta|u|^pdx\, ,
\nonumber
\end{eqnarray}
for all $u\in C^{\infty}_c(\Omega\setminus K)$; here $\zeta=\zeta(d(x)/D)$.
\label{petr}
\end{prop}
{\em Proof.} Let $u\in C^{\infty}_c(\Omega\setminus K)$ be fixed.
For a positive, locally bounded function $\phi$ with $|\nabla\phi|\in L^1_{\rm loc}(\Omega\setminus K)$ we have
\begin{eqnarray*}
-\int_{\Omega}\Delta\phi |u|^pdx&=&p\int_{\Omega}\nabla\phi\cdot
(|u|^{p-2}u\nabla u)dx \\
&=&-p\int_{\Omega}\phi |u|^{p-2}u\Delta u\, dx -p(p-1)\int_{\Omega}\phi |u|^{p-2}|\nabla u|^2 dx \\
&\leq& p\Big (\frac{p-1}{p}\int_{\Omega}d^{\frac{\gamma}{p-1}}(1+\mu\zeta)^{-\frac{1}{p-1}}\phi^{\frac{p}{p-1}}|u|^pdx + \\
&& \quad  +\frac{1}{p}\int_{\Omega}d^{-\gamma}(1+\mu\zeta)|\Delta u|^pdx\Big)
 -p(p-1)\int_{\Omega}\phi |u|^{p-2}|\nabla u|^2 dx, 
\end{eqnarray*}
from which follows that
\begin{equation}
 \int_{\Omega}d^{-\gamma}(1+\mu\zeta) |\Delta u|^pdx \geq  T_1+T_2+T_3 \; ,
\label{pasxa}
\end{equation}
where
\begin{eqnarray*}
&& T_1=p(p-1)\int_{\Omega}\phi |u|^{p-2}|\nabla u|^2dx\, , \\
&& T_2=- \int_{\Omega}\Delta\phi |u|^p dx\, ,\\
&& T_3 =-(p-1)\int_{\Omega}d^{\frac{\gamma}{p-1}}(1+\mu\zeta)^{-\frac{1}{p-1}}\phi^{\frac{p}{p-1}}|u|^pdx\, .
\end{eqnarray*}
We next choose $\phi=\lambda d^{-\gamma -2p+2}(1+\alpha\eta+\beta\eta^2)$
where $\lambda>0$ and $\alpha,\beta\in\R$
are to be determined and $\eta=\eta(d(x)/D)$ with $D\geq \sup_{\Omega}d$ also yet to be determined. 
To estimate $T_1$ we set $v=|u|^{p/2}$ and apply \cite[Theorem 1]{BT} obtaining
\begin{eqnarray}
T_1&=&\frac{4(p-1)\lambda}{p}\int_{\Omega}d^{-\gamma-2p+2}(1+\alpha\eta+\beta\eta^2)
|\nabla v|^2dx \nonumber\\
&\geq& \frac{4(p-1)\lambda}{p}\int_{\Omega}d^{-\gamma -2p}
\bigg\{\frac{(k-\gamma -2p)^2}{4}
+\frac{(k-\gamma -2p)^2\alpha}{4}\eta + \nonumber \\
&& \qquad + \biggl(\frac{(k-\gamma -2p)\alpha}{4}+\frac{(k-\gamma -2p)^2\beta}{4}\biggr)
\eta^2 + \label{t1} \\
&&\qquad + \biggl(\frac{1}{4}+
\frac{(k-\gamma -2p)\alpha}{4}\biggr)\zeta \bigg\}|u|^pdx\, .
\nonumber
\end{eqnarray}
To estimate $T_2$ we define $f(t)=
\lambda t^{-\gamma -2p+2}(1+\alpha\eta(t/D)+\beta\eta^2(t/D))$, $t\in (0,\sup_{\Omega}d)$,
so that $\phi(x)=f(d(x))$. We then have
\begin{eqnarray}
&&f'(t) =\lambda t^{-\gamma-2p+1}\Big\{ (-\gamma-2p+2)
+  (-\gamma-2p+2)\alpha \eta +\Big[\frac{\alpha}{2}+(-\gamma-2p+2)\beta\Big]\eta^2 \nonumber \\
&&\qquad\qquad\qquad \qquad  +\frac{\alpha}{2}\zeta 
+\beta\eta^3 +\beta\eta\zeta\Big\}
\label{gio}
\end{eqnarray}
and
\begin{eqnarray}
f''(t) &=&\lambda t^{-\gamma-2p}\bigg\{ (\gamma+2p-1)(\gamma+2p-2)
+(\gamma+2p-1)(\gamma+2p-2)\alpha\eta \nonumber\\
&&+\Big[  -\frac{(2\gamma+ 4p -3)\alpha}{2} + (\gamma+2p-1)(\gamma+2p-2)\beta \Big]\eta^2 \label{kou} \\
&&+\Big[  -\frac{(2\gamma+ 4p -3)\alpha}{2} \Big]\zeta +
\Big[ \frac{\alpha}{2} -(2\gamma+ 4p -3)\beta\Big](\eta^3+\eta\zeta) + \alpha\theta\bigg\}
\nonumber
\end{eqnarray}
(where the argument of $\eta,\zeta$ and $\theta$ in (\ref{gio}) and (\ref{kou}) is $t/D$).
Since $f'(t)\leq 0$ for large $D$, we have from (\ref{woody})
\begin{equation}
-\Delta\phi =-f''(d)- f'(d)\Delta  d \geq -f''(d)-\frac{k-1}{d}f'(d) \; , 
\label{tel}
\end{equation}
in the distributional sense in $\Omega\setminus K$.
Combining (\ref{gio}), (\ref{kou}) and (\ref{tel}) we conclude that
\begin{eqnarray}
T_2&\geq&\lambda \int_{\Omega}d^{-\gamma-2p}\Biggl\{ (k-\gamma-2p)(\gamma +2p-2)+
(k-\gamma-2p)(\gamma +2p-2) \alpha\eta +\nonumber\\
&&\qquad+\Big( \frac{2\gamma +4p-k-2}{2}\alpha + (k-\gamma-2p)(\gamma +2p-2)\beta  \Big)\eta^2 \nonumber\\
&&\qquad+\Big( \frac{2\gamma +4p-k-2}{2}\alpha\Big)\zeta  \label{t2}\\
&&\qquad +\Big( (2\gamma +4p-k-2)\beta-\frac{\alpha}{2}\Big)(\eta^3+\eta\zeta) - \alpha\theta +O(\eta^4)\Biggr\}dx.
\nonumber
\end{eqnarray}
As for $T_3$, we use Taylor's theorem to obtain after some simple calculations
\begin{eqnarray*}
&&(1+\mu\zeta)^{-\frac{1}{p-1}}\phi^{\frac{p}{p-1}}\\
&=&\lambda^{\frac{p}{p-1}}
d^{-\frac{(\gamma+2p-2)p}{p-1}}\bigg\{
1+\frac{p\alpha}{p-1}\eta +
\biggl(\frac{p\beta}{p-1}+\frac{p\alpha^2}{2(p-1)^2}\biggr)\eta^2 +\\
&& -\frac{\mu}{p-1}\zeta +\biggl(\frac{p\alpha\beta}{(p-1)^2}-\frac{p(p-2)\alpha^3}{6(p-1)^3}
\biggr)\eta^3  -\frac{p\alpha\mu}{(p-1)^2}\eta\zeta  + O(\eta^4)\bigg\},
\end{eqnarray*}
and thus conclude that
\begin{eqnarray}
T_3&=&\!\!-(p-1)\lambda^{\frac{p}{p-1}}\int_{\Omega}d^{-\gamma -2p}
\Biggl\{1+\frac{p\alpha}{p-1}\eta +
\biggl(\frac{p\beta}{p-1}+\frac{p\alpha^2}{2(p-1)^2}\biggr)\eta^2 -\frac{\mu}{p-1}\zeta + \nonumber\\
&&\hspace{2cm}+\biggl(\frac{p\alpha\beta}{(p-1)^2}-\frac{p(p-2)\alpha^3}{6(p-1)^3}
\biggr)\eta^3 - \frac{p\alpha\mu}{(p-1)^2} \eta\zeta  +O(\eta^4)\Biggr\}dx.
\label{t3}
\end{eqnarray}
Using (\ref{t1}), (\ref{t2}) and (\ref{t3}) we arrive at
\begin{equation}
\int_{\Omega}d^{-\gamma}(1+\mu\zeta)|\Delta u|^pdx\geq \int_{\Omega}d^{-\gamma-2p}V|u|^pdx
\label{eq:v}
\end{equation}
where the function $V$ has the form
\[
V=r_0+r_1\eta +r_2\eta^2 +r_2'\zeta +r_3\eta^3+r_3'\eta\zeta+r_3''\theta +O(\eta^4).
\]
We compute the coefficients $r_i,r_i',r_i''$ by adding the respective coefficients from
(\ref{t1}), (\ref{t2}) and (\ref{t3}). We find
\begin{eqnarray}
&&r_0=\frac{(k-\gamma -2p)(pk-k+\gamma)}{p}\lambda-(p-1)\lambda^{\frac{p}{p-1}}\; , \nonumber\\
&&r_1=\frac{(k-\gamma -2p)(pk-k+\gamma)}{p}\lambda\alpha -p\lambda^{\frac{p}{p-1}}\alpha \; , \nonumber\\
&&r_2= \frac{pk+2p-2k+2\gamma}{2p}\alpha\lambda
+\frac{(k-\gamma-2p)(pk-k+\gamma)}{p}\beta\lambda -\nonumber\\
&&\qquad\qquad\qquad\quad -(p-1)\biggl(\frac{p\beta}{p-1}+
\frac{p\alpha^2}{2(p-1)^2}\biggr)\lambda^{\frac{p}{p-1}} \; , \label{ooo}\\
&&r_2'=\biggl(\frac{p-1}{p}+\frac{pk-2k+2p+2\gamma}{2p}\alpha\biggr)\lambda +\mu\lambda^{\frac{p}{p-1}} \; , \nonumber\\
&&r_3=\lambda \Big(-\frac{\alpha}{2}+(2\gamma+4p-k-2)\beta\Big) -(p-1)\lambda^{\frac{p}{p-1}}
\Big( \frac{p\alpha\beta}{(p-1)^2}-\frac{p(p-2)\alpha^3}{6(p-1)^3}\Big) \; , \nonumber\\
&&r_3'=\lambda \Big(-\frac{\alpha}{2}+(2\gamma+4p-k-2)\beta\Big) 
+\frac{p\alpha\mu}{p-1}\lambda^{\frac{p}{p-1}} \; , \nonumber\\
&&r_3''=-\lambda\alpha\; .
\nonumber
\end{eqnarray}
We now proceed to specify $\lambda$, $\alpha$ and $\beta$. We choose
$\lambda$ so as to optimize $r_0$, which yields
\[
\lambda=Q^{p-1}\quad , \quad
r_0 =Q^p.
\]
Then $r_1=0$ irrespective of the choice of $\alpha$ and $\beta$.
We subsequently choose
\[
\alpha =\frac{(p-1)(pk-2k+2p+2\gamma)}{(k-\gamma-2p)(pk-k+\gamma)},
\]
which yields $r_2=0$ and $r_2'=B(2,\gamma)+\mu A(2,\gamma)$. Hence it remains to show that
$\beta$ can be chosen so that for large enough $D$ there holds
\begin{equation}
r_3\eta^3+r_3'\eta\zeta+r_3''\theta +O(\eta^4)\geq 0\quad  \mbox{ in }\Omega.
\label{sss}
\end{equation}
This is done in the following lemma and this is where condition $\tc$ is needed. $\hfill //$

\begin{lemma}
If $\gamma\neq (3pk-8p^2-2k+6p)/(4p-2)$ or $p> (13+\sqrt{105})/4$
then there exists $\beta\in\R$ such that for large enough $D$ there holds
\begin{equation}
r_3\eta^3+r_3'\eta\zeta+r_3''\theta +O(\eta^4)\geq 0\quad  \mbox{ in }\Omega
\label{sss1}
\end{equation}
(here $\eta=\eta(d(x)/D)$, and similarly for $\zeta$ and $\theta$).
\label{lem:new}
\end{lemma}
{\em Proof.} We claim that it is enough to find $\beta\in\R$ such that for large enough $D$ we have
\begin{equation}
r_3+r_3'+r_3''>0\; .
\label{opo}
\end{equation}
Indeed, the fact that
\[
\lim_{t\to 0+}\frac{\eta^3(t)}{X_1^3(t)}=\lim_{t\to 0+}\frac{\eta(t)\zeta(t)}{X_1^3(t)}
=\lim_{t\to 0+}\frac{\theta(t)}{X_1^3(t)}=1\; , 
\]
implies that
\[
r_3\eta^3+r_3'\eta\zeta+r_3''\theta +O(\eta^4)=\Big(r_3+r_3'(1+o(1)) +r_3''(1+o(1))\Big)\eta^3+O(\eta^4) 
\]
where $\lim o(1)=0$ as $D\to +\infty$, uniformly in $x\in\Omega$; hence (\ref{sss1}) follows.

To prove (\ref{opo}) we calculate $r_3,r_3'$ and $r_3''$; from (\ref{ooo}) we obtain
\begin{eqnarray*}
&&r_3=\Big( -\frac{\alpha Q^{p-1}}{2}+\frac{p(p-2)\alpha^3Q^p}{6(p-1)^2}\Big) - R \beta\; , \\
&&r_3'=\Big( -\frac{\alpha Q^{p-1}}{2}+\frac{p\mu\alpha Q^p}{p-1}\Big) 
 + (-R +\frac{p\alpha Q^p}{p-1} )\beta \; , \\
&&r_3''=-\alpha Q^{p-1} \; ,
\end{eqnarray*}
where $R=2(p-1)(k-\gamma-2p)Q^{p-1}/p\, .$
We distinguish two cases. \nl
$\ia$ $\gamma\neq (3pk-8p^2-2k+6p)/(4p-2)$.
We then observe that the coefficient of $\beta$ in
$r_3+r_3'+r_3''$ is non-zero. Hence (\ref{opo}) is satisfied if $\beta$ is either large and negative
or large and positive.

$\ib$ $\gamma = (3pk-8p^2-2k+6p)/(4p-2)$. We then choose $\beta=0$ and we have
\[
Q=\frac{(k-2)^2(4p-3)}{4(2p-1)^2}\quad , \quad \alpha Q=\frac{2(k-2)(p-1)^2}{p(2p-1)},
\]
from which follows that
\[
r_3+r_3'+r_3''=\frac{2(2p-1)\alpha Q^{p-1}}{3p(4p-3)}(2p^2-13p+8)\; .
\]
Since $\alpha>0$ in this case, this is positive as $(13+\sqrt{105})/4$ is the largest root
of the polynomial $2p^2-13p+8$. $\hfill //$

{\em Note.}
Using $\phi=\lambda d^{-\gamma -2p+2}(1+\alpha\eta+\beta\eta^2+\beta_1\zeta)$ in order to remove
$\tc$ does not work, as the coefficient of $\beta_1$ in $r_3+r_3'+r_3''$ turns out to be zero when the corresponding coefficient
of $\beta$ is zero.

\begin{lemma}
Let $m\in\N$ and $\gamma\geq 0$. Then:
\begin{eqnarray*}
&\ia& \mbox{ If $m$ is even then }\\
&& (a)\quad A(m,\gamma)=A(2,\gamma)A(m-2,\gamma+2p)  , \\
&& (b) \quad  B(m,\gamma)=A(2,\gamma)B(m-2,\gamma+2p) +A(m-2,\gamma+2p)B(2,\gamma) . \\[0.2cm]
&\ib& \mbox{ If $m$ is odd then }\\
&& (a)\quad  A(m,\gamma)=A(1,\gamma)A(m-1,\gamma+p)   , \\
&& (b)\quad B(m,\gamma)=A(1,\gamma)B(m-1,\gamma+p) +A(m-1,\gamma+p)B(1,\gamma).
\end{eqnarray*}
\label{lem:papa}
\end{lemma}
{\em Proof.} We shall only prove $\ia$(b), the other cases being simpler or similar. 
So let us assume that $m=2r$, $r\in\N$. Then
\begin{eqnarray*}
&& \hspace{-1cm}A(2,\gamma)B(2r-2,\gamma+2p)+A(2r-2,\gamma+2p)B(2,\gamma)\\
&=&\Big(\frac{k-\gamma-2p}{p}\Big)^p \Big(\frac{pk-k+\gamma}{p}\Big)^p \; \frac{p-1}{2p}\\
&&\qquad \times\prod_{i=0}^{r-2}\prod_{j=1}^{r-1}\Big(\frac{k-\gamma-(2r-2i)p}{p}\Big)^p
\Big(\frac{pk-k+\gamma+(2r-2j)p}{p}\Big)^p\\
&&\qquad \times\bigg\{ \sum_{i=0}^{r-2}\Big(\frac{k-\gamma-(2r-2i)p}{p}\Big)^{-2} +
\sum_{j=1}^{r-1}\Big(\frac{pk-k+\gamma+(2r-2j)p}{p}\Big)^{-2}\bigg\}\\
&&\; +\prod_{i=0}^{r-2}\prod_{j=1}^{r-1}\Big(\frac{k-\gamma-(2r-2i)p}{p}\Big)^p
\Big(\frac{pk-k+\gamma+(2r-2j)p}{p}\Big)^p\\
&&\qquad \times \frac{p-1}{2p}\Big(\frac{k-\gamma-2p}{p}\Big)^p \Big(\frac{pk-k+\gamma}{p}\Big)^p 
\bigg\{ \Big(\frac{k-\gamma-2p}{p}\Big)^{-2} + \Big(\frac{pk-k+\gamma}{p}\Big)^{-2} \bigg\}\\
&=&\frac{p-1}{2p}\prod_{i=0}^{r-1}\prod_{j=1}^{r}\Big(\frac{k-\gamma-(2r-2i)p}{p}\Big)^p
\Big(\frac{pk-k+\gamma+(2r-2j)p}{p}\Big)^p \\
&&\qquad \times \bigg\{ \sum_{i=0}^{r-1}\Big(\frac{k-\gamma-(2r-2i)p}{p}\Big)^{-2} +
\sum_{j=1}^{r}\Big(\frac{pk-k+\gamma+(2r-2j)p}{p}\Big)^{-2}\bigg\}\\
&=&A(2r,\gamma),
\end{eqnarray*}
as claimed. $\hfill //$

{\bf Proof of Theorem \ref{thmb}}
Before proceeding with the proof, let us make a comment on its assumptions.
The proof essentially uses iteration. For example, if $m$ is even, then we repeatedly use
Proposition \ref{petr} obtaining
\[
\int_{\Omega}\frac{|\Delta^{m/2}u|^p}{d^{\gamma}}dx \geq \int_{\Omega}(a_1+b_1\zeta) \frac{|\Delta^{(m-2)/2}u|^p}{d^{\gamma+2p}}dx \geq
\int_{\Omega}(a_2+b_2\zeta) \frac{|\Delta^{(m-4)/2}u|^p}{d^{\gamma+4p}}dx \geq \ldots\, ,
\]
etc. Hence at the $(i+1)$th step, $0\leq i\leq (m-2)/2$, we estimate the integral
$\int_{\Omega}(a_i+b_i\zeta)  d^{-(\gamma+2ip)}|\Delta^{(m-2i)/2}u|^pdx$.
In applying Proposition \ref{petr}, we verify that $\ia$ $k-(\gamma+2ip)-2p>0$
(this is satisfied since $k-\gamma-mp>0$) and $\ib$ if $p\leq  (13+\sqrt{105})/4$,
then $\gamma +2ip\neq (3pk-8p^2-2k+6p)/(4p-2)$. This is indeed the case by the assumption of the theorem since
$\gamma +jp >3pk-8p^2-2k+6p)/(4p-2)$ for any $j\geq 2$ (recall that $\gamma\geq 0$).

We now come to the details of the proof.
We shall use induction on $[(m+1)/2]$. If $[(m+1)/2]=1$, that is $m=1$ or $m=2$,
then (\ref{1.2}) follows from \cite[Theorem 1]{BT} or Proposition \ref{petr} respectively. We assume
that the statement of the theorem is valid for $[(m+1)/2]\in \{1,2,\ldots,r-1\}$ and
consider the case $[(m+1)/2]=r$. 
For this we distinguish two cases, depending on whether $m$ is even or odd.

$\ia$ $m$ even. We first use Proposition \ref{petr} and then the induction hypothesis
(and for this we note that the assumption $k-\gamma-mp>0$ implies both
$k-\gamma-2p>0$ and $k-(\gamma+2p)-(m-2)p>0$). We have
\begin{eqnarray*}
&&\int_{\Omega}d^{-\gamma}(1+\mu\zeta)|\Delta^{m/2}u|^pdx \\
&\geq&A(2,\gamma)
\int_{\Omega}d^{-\gamma-2p}|\Delta^{\frac{m-2}{2}}u|^pdx
+ [B(2,\gamma) + A(2,\gamma)\mu ]\int_{\Omega}d^{-\gamma-2p}\zeta|\Delta^{\frac{m-2}{2}}u|^pdx \\
&\geq&\!A(2,\gamma)\bigg\{ A(m-2,\gamma+2p)\int_{\Omega}\! d^{-\gamma-mp}|u|^pdx
+\! B(m-2,\gamma+2p)\! \int_{\Omega}\! d^{-\gamma-mp}\zeta |u|^pdx\bigg\}\\
&&+[B(2,\gamma) + A(2,\gamma)\mu ]A(m-2,\gamma+2p)\int_{\Omega}d^{-\gamma-mp}\zeta|u|^pdx\\
&=& A(2,\gamma)A(m-2,\gamma+2p) \int_{\Omega}d^{-\gamma-mp}|u|^pdx +\\
&& + \bigg\{ \Big[A(2,\gamma)B(m-2,\gamma+2p) +A(m-2,\gamma+2p)B(2,\gamma)\Big] + \\
&&+ A(2,\gamma)A(m-2,\gamma+2p)\mu \bigg\} \int_{\Omega}d^{-\gamma-mp}\zeta |u|^pdx,
\end{eqnarray*}
and the proof is complete if we recall Lemma \ref{lem:papa}.\nl
$\ib$ $m$ odd. The proof is similar, the only difference being that we use \cite[Theorem 1]{BT} instead of Proposition \ref{petr}.
We omit the details. $\hfill //$

\begin{remark}
{\rm We point out that in the proofs of Proposition \ref{petr} and Theorem \ref{thmb}
we did not use at any point the assumption that $k$ is the codimension of the set $K$. Indeed, a
careful look at the two proofs shows that $K$ can be any closed
set such that $\dist(x,K)$ is bounded in $\Omega$ and for which
the inequality  $d\Delta d -k+1\geq 0$ is satisfied in $\Omega\setminus K$; the proof does not even require $k$ to
be an integer. Of course, the natural realization of this
assumption is that $K$ is smooth and $k=\codim(K)$.}
\end{remark}

Let us define the {\em inradius of $\Omega$ relative to $K$} by ${\rm Inr}(\Omega ; K)=\sup_{\Omega}d(x)$.
Looking at the proof of Theorem \ref{thmb} we see that when $D$ is chosen large enough, the actual requirement 
is that $d(x)/D$ is small uniformly in $x\in\Omega$. This, combined with the fact that $t^{-\gamma-mp}X_1^2\ldots X_i^2(t)$
has a positive minimum in $(0,1)$, leads to the following corollary of Theorem \ref{thmb}:

\begin{coro}
Under the conditions of Theorem \ref{thmb} for any $r\geq 0$ there exists a constant $c=c(m,p,k,r)>0$ such that
\begin{eqnarray}
&&\int_{\Omega}d^{-\gamma} |\Delta^{m/2}u|^pdx \geq
A(m,\gamma)\int_{\Omega}d^{-\gamma-mp}|u|^p dx +\label{1.9}\\[0.2cm]
&&\qquad +B(m,\gamma)\sum_{i=1}^r\int_{\Omega}d^{-\gamma -mp}
 X_1^2X_2^2\ldots X_i^2 |u|^p dx +c\, {\rm Inr}(\Omega ; K)^{-\gamma -mp}\int_{\Omega}|u|^pdx,
\nonumber
\end{eqnarray}
for all $u\in C^{\infty}_c(\Omega\setminus K)$.
\end{coro}

We end this section with a proposition about the case where condition $\tc$ is not satisfied.

\begin{prop}
Suppose that all conditions of Theorem \ref{thmb} except $\tc$ are satisfied. Then
\begin{equation}
\int_{\Omega}d^{-\gamma} |\Delta^{m/2}u|^pdx \geq
A(m,\gamma)\int_{\Omega}d^{-\gamma-mp}|u|^p dx + c_{\epsilon} \int_{\Omega}d^{-\gamma-mp+\epsilon}|u|^p dx\; ,
\label{new}
\end{equation}
for any $\epsilon>0$ and all $u\in C^{\infty}_c(\Omega\setminus K)$.
\end{prop}
{\em Proof.} We only give a sketch of the proof. Suppose first that $m=2$.
We use (\ref{pasxa}), but this time with $\phi=\lambda d^{-\gamma-2p+2}(1+\mu d^{\epsilon})$;
here $\mu$ is to be determined and $\lambda=Q^{p-1}$. Arguing as in the proof of Proposition \ref{petr} we
obtain
\begin{equation}
\int_{\Omega}d^{-\gamma} |\Delta^{m/2}u|^pdx \geq
A(m,\gamma)\int_{\Omega}d^{-\gamma-mp}|u|^p dx + \int_{\Omega}\tilde{V}d^{-\gamma-2p}|u|^pdx\; ,
\label{new1}
\end{equation}
where
\begin{eqnarray*}
\tilde{V}(x)&=&\frac{\lambda\mu}{p}(p-1)(k-\gamma-2p+\epsilon)^2d^{\epsilon} +\lambda\mu(k-\gamma-2p+\epsilon)
(\gamma+2p-2-\epsilon)d^{\epsilon}\\
&&+(p-1)\lambda^{\frac{p}{p-1}} -(p-1)\lambda^{\frac{p}{p-1}}(1+\mu d^{\epsilon})^{\frac{p}{p-1}};
\end{eqnarray*}
here we have added and subtracted $(p-1)\lambda^{\frac{p}{p-1}}$ in order to create the first term in the
right-hand side of (\ref{new1}).
Using Taylor's
theorem we obtain
\[
\tilde{V}(x)=(c_1\lambda\mu\epsilon +O(\epsilon^2))d^{\epsilon} +O(d^{2\epsilon}),
\]
where $c_1=(pk-2k+2\gamma+2p)/p$. The fact that $\tc$
is violated implies that $c_1\neq 0$, and choosing $\mu$ so that $c_1\mu>0$ completes
the proof when $m=2$. Iteration yields the result in the general case when $m$ is even.
The case where $m$ is odd is treated similarly. $\hfill //$

\section{Optimality of the constants} 

In this section we present the proof of Theorem \ref{thm:meli}. Hence we assume throughout that
$\Omega$ is domain in $\R^N$ and $K$ is an affine hypersurface of codimension $k\in\{1,2,\ldots,N\}$.

For the sake of simplicity we shall only consider the special case $\gamma=0$, the proof in the general case
presenting no difference whatsoever other than the additional dependence of some constants on $\gamma$.
Also, for the sake of bravity we shall prove the theorem only for $m$ even, the proof when $m$ is odd being similar.

Hence, writing $A(2m)$ and $B(2m)$ for $A(2m,0)$ and
$B(2m,0)$ respectively, we intend to look closely at
\[
I_{2m,r-1}[u]:=\int_{\Omega}|\Delta^mu|^pdx -|A(2m)|\int_{\Omega}\frac{|u|^p}{d^{2mp}}dx -|B(2m)|
\sum_{i=1}^{r-1}\int_{\Omega}\frac{|u|^p}{d^{2mp}}X_1^2\ldots X_i^2dx,
\]
for particular test functions $u$; here and below, $X_j=X_j(d(x)/D)$ for some fixed $D\geq\sup_{\Omega}d(x)$.
We begin by defining the polynomial
\[
\alpha_m(s)= \prod_{i=0}^{m-1} (s-2i) \prod_{j=1}^{m}(s+k-2j) \;, \quad s\in\R\, ,
\]
which will play an important role in the sequel.
\begin{lemma}
There holds
\begin{eqnarray}
\ia &&  |A(2m)|=|\alpha_m|^p\bigg|_{s=\frac{2mp-k}{p}} \; , \label{av1}\\
\ib &&  |B(2m)|= \frac{p-1}{2p}|\alpha_m|^{p-2}(\alpha_m'^2-\alpha_m\alpha_m'')\bigg|_{s=\frac{2mp-k}{p}}.\label{av2}
\end{eqnarray}
\label{radio}
\end{lemma}
{\em Proof.} Part $\ia$ is easily verified. From the relation
\[
\frac{\alpha_m'}{\alpha_m}= \sum_{i=0}^{m-1}(s-2i)^{-1} +\sum_{j=1}^m(s+k-2j)^{-1},
\]
we obtain
\[
\alpha_m^{-2}(\alpha_m'^2-\alpha_m\alpha_m'')
=-\Big( \frac{\alpha_m'}{\alpha_m}\Big)'= \sum_{i=0}^{m-1}(s-2i)^{-2} +\sum_{j=1}^m(s+k-2j)^{-2},
\]
and $\ib$ follows. $\hfill //$

Let $s_0>(2mp-k)/p$ and $s_1,\ldots,s_r\in\R$ be fixed parameters. For $0\leq i\leq j\leq r$ we define
\[
Y_{ij}=X_1^2\ldots X_i^2 X_{i+1}\ldots X_j \; ,
\]
with the natural interpretations $Y_{00}=1$, $Y_{ii}=X_1^2\ldots X_i^2$, $Y_{0j}=X_1\ldots X_j$. We then define the integrals
\begin{eqnarray*}
\Gamma_{ij}&=&\int_{\Omega}d^{(s_0-2m)p}X_1^{ps_1}\ldots X_r^{ps_r}Y_{ij}dx \\
&=&\int_{\Omega}d^{(s_0-2m)p} X_1^{ps_1+2}\!\!\!\!\ldots
X_i^{ps_i+2}X_{i+1}^{ps_{i+1}+1}\!\!\ldots X_j^{ps_j+1}
X_{j+1}^{ps_{j+1}}\!\!\!\!\ldots X_r^{ps_r} dx.
\end{eqnarray*}
\begin{lemma}
Let $u(x)=d^{s_0}X_1^{s_1}\ldots X_r^{s_r}$, where $X_i=X_i(d(x)/D)$. Then
\begin{equation}
I_{2m,r-1}[u]=\sum_{0\leq i\leq j\leq r}a_{ij}\Gamma_{ij} + 
\int_{\Omega}d^{(s_0-2m)p}X_1^{ps_1}X_2^{ps_2}\ldots X_r^{ps_r}O(X_1^3)dx \; ,
\label{rain}
\end{equation}
where
\begin{eqnarray}
a_{00} &=&|\alpha_m|^p -|A(2m)| , \nonumber\\
a_{0j} &=& ps_j|\alpha_m|^{p-2}\alpha_m\alpha_m' ,  \quad  1\leq j\leq r,  \nonumber\\
a_{ii} &=&  \frac{ps_i}{2}|\alpha_m|^{p-2}\Big( \alpha_m\alpha_m''(s_i+1) +(p-1)\alpha_m'^2 s_i\Big) -|B(2m)| , \quad   1\leq i\leq r-1, \nonumber \\
a_{rr}&=& \frac{ps_r}{2}|\alpha_m|^{p-2}\Big( \alpha_m\alpha_m''(s_r+1) +(p-1)\alpha_m'^2s_r\Big) ,\label{aij}\\
a_{ij}& =& \frac{ps_j}{2}|\alpha_m|^{p-2}\Big( \alpha_m\alpha_m''(2s_i+1) +2(p-1)\alpha_m'^2s_i\Big) , 
\quad 1\leq i<j\leq r\; ; \nonumber
\end{eqnarray}
here and below, $\alpha_m$, $\alpha_m'$ and $\alpha_m''$ stand for $\alpha_m(s_0)$, $\alpha_m'(s_0)$ and $\alpha_m''(s_0)$
respectively.
\label{trian}
\end{lemma}
{\em Proof.} The fact that $K$ is affine implies that $\Delta d=(k-1)/d$ and therefore
\begin{equation}
\Delta (f(d))=f''(d)+\frac{k-1}{d}f'(d)\; ,
\label{ind}
\end{equation}
for any smooth function $f$ on $(0,+\infty)$.
We define the functions $g$, $\tilde{g}$ by
\[
g(x)=s_1X_1+s_2X_1X_2 +\ldots s_rX_1X_2\ldots X_r \; , \quad \nabla g=\frac{\tilde{g}}{d}\nabla d\; , 
\]
and observe that by (\ref{diffrule}),
\begin{equation}
g^3(t)=O(X_1^3) \quad , \quad \tilde{g}^2(t)=O(X_1^4).
\label{snow}
\end{equation}
Now, (\ref{ind}) and (\ref{snow}) together with a simple induction argument on $m$ imply
\[
\Delta^m u =d^{s_0-2m}X_1^{s_1}\ldots X_r^{s_r}\Big(\alpha_m + \alpha_m' g(d) +\frac{\alpha_m''}{2}g^2(d) +
\frac{\alpha_m''}{2}\tilde{g}(d) +O(X_1^3)\Big) .
\]
Using Taylor's theorem we then obtain
\begin{eqnarray}
&& |\Delta^m u|^p =d^{(s_0-2m)p}X_1^{ps_1}\ldots X_r^{ps_r}\Bigg\{|\alpha_m|^p + p|\alpha_m|^{p-2}\alpha_m\alpha_m' g(d)
\label{dr}\\
&&\qquad\qquad \quad +\frac{p}{2}|\alpha_m|^{p-2}\Big( \alpha_m\alpha_m'' +(p-1)\alpha_m'^2\Big)g^2
+\frac{p}{2}|\alpha_m|^{p-2}\alpha_m\alpha_m''\tilde{g} +O(X_1^3)\Bigg\}.
\nonumber
\end{eqnarray}
On the other hand we have  (cf. (\ref{diffrule}))
\begin{eqnarray}
&& \int_{\Omega}d^{(s_0-2m)p}X_1^{ps_1}\ldots X_r^{ps_r}g\, dx =\sum_{j=1}^rs_j\Gamma_{0j}, \nonumber\\
&& \int_{\Omega}d^{(s_0-2m)p}X_1^{ps_1}\ldots X_r^{ps_r}g^2 dx=\sum_{i=1}^rs_i^2\Gamma_{ii}
 +2\sum_{1\leq i<j\leq r}s_is_j\Gamma_{ij}, \label{net}\\
&& \int_{\Omega}d^{(s_0-2m)p}X_1^{ps_1}\ldots X_r^{ps_r}\tilde{g}\, dx
=\sum_{i=1}^rs_i\Gamma_{ii} +\sum_{1\leq i<j\leq r}s_j\Gamma_{ij}\; .
\nonumber
\end{eqnarray}
The stated relation follows from (\ref{dr}),  (\ref{net}) and the fact that
$I_{2m,r-1}[u]=\int_{\Omega}|\Delta^m u|^pdx -|A(2m)|\Gamma_{00}-|B(2m)|\sum_{i=1}^{r-1}\Gamma_{ii}$.
$\hfill //$

Up to this point the exponents $s_0,s_1,\ldots, s_r$ where arbitrary subject only to
$s_0>(2mp-k)/p$. We now make a more specific choice, taking
\begin{equation}
s_0=\frac{2mp-k+\epsilon_0}{p}\; , \qquad s_j=\frac{-1+\epsilon_j}{p}\; , \quad 1\leq j\leq r\, ,
\label{se}
\end{equation}
where $\epsilon_0,\ldots,\epsilon_r$ are small positive parameters.
We consider $I_{2m,r-1}[u]$ as a function of these parameters and intend
to take the limits $\epsilon_0\searrow 0,\ldots ,\epsilon_r\searrow 0$. In taking these
limits we shall ignore terms that are bounded uniformly in the $\epsilon_i$'s. In order to distinguish
such terms we shall need the following criterion, which is a simple consequence of (\ref{diffrule}):
\begin{equation}
\int_{\Omega} d^{-k+\epsilon_0}X_1^{1+\epsilon_1}\ldots X_r^{1+\epsilon_r}dx <\infty 
\Longleftrightarrow\; \left\{
\begin{array}{ll}
& \epsilon_0>0 \\
\mbox{ or}&\mbox{$\epsilon_0=0$ and $\epsilon_1>0$}\\
\mbox{ or}&\mbox{$\epsilon_0=\epsilon_1=0$ and $\epsilon_2>0$}\\
& \cdots \\
\mbox{ or}&\mbox{$\epsilon_0=\epsilon_1=\ldots =\epsilon_{r-1}=0$ and $\epsilon_r>0$.}
\end{array}\right.
\label{finite}
\end{equation}
Also, concerning terms that diverge as the $\epsilon_i$'s tend to zero,
we shall need some quantitive information on the rate of divergence as well as some
mutual cancelation properties. These are collected in the following
\begin{lemma}
We have
\begin{eqnarray*}
&\ia& \int_{\Omega}d^{-k+\epsilon_0}X_1^{\beta}dx
\leq c_{\beta}\epsilon_0^{-1+\beta}\; ,\quad \beta <1 ;\\
&\ib & \int_{\Omega}d^{-k}X_1\ldots X_{i-1}X_i^{1+\epsilon_i}X_{i+1}^{\beta}dx
\leq c_{\beta}\epsilon_i^{-1+\beta}\; , \qquad \beta <1 \; , \quad 1\leq i\leq r-1\, ; \\
&\ic &  \epsilon_0^2\Gamma_{00} -2\epsilon_0\sum_{j=i+1}^r(1-\epsilon_j)\Gamma_{0j} =\\
&&\qquad \qquad
\sum_{i=1}^r(\epsilon_i-\epsilon_i^2)\Gamma_{ii} -\sum_{1\leq i<j\leq r}(1-\epsilon_j)(1-2\epsilon_i)\Gamma_{ij}
+O(1), \\
&&\mbox{where the O(1) is uniform in }\epsilon_0,\ldots,\epsilon_r \; ;\\
& \id & \mbox{let $i\geq 0$ and (if $i\geq 1$) assume that $\epsilon_0=\ldots=\epsilon_{i-1}=0$. Then}\\
&&\qquad\qquad\qquad \epsilon_i\Gamma_{ii} = \sum_{j=i+1}^r(1-\epsilon_j)\Gamma_{ij} +O(1), \\
&&\mbox{where the O(1) is uniform in $\epsilon_i,\ldots,\epsilon_r$.}
\end{eqnarray*}
\label{diverge}
\end{lemma}
{\em Proof.} Parts $\ia$ and $\ib$ are proved using the coarea formula and \cite[Lemma 9]{B}.
Parts $\ic$ and $\id$ are proved by integrating by parts; see \cite{BFT2}, pages 181 and 184 respectively
for the detailed proof. $\hfill//$

\begin{remark}
{\rm We are now in position to prove Theorem \ref{thm:meli}, but before proceeding
some comments are necessary. The proof of the theorem is local: we fix a point $x_0\in \Omega\cap K$ and work
entirely in a small ball $B(x_0,\delta)$ using a cut-off function $\phi$. The sequence of functions
that is used is then given by
\[
u(x)=\phi(x) d(x)^{\frac{2mp-k+\epsilon_0}{p}}X_1(d(x)/D)^{\frac{-1+\epsilon_1}{p}}\ldots
X_r(d(x)/D)^{\frac{-1+\epsilon_r}{p}}, \qquad (\epsilon_0,\ldots,\epsilon_r>0)
\]
and, as already mentioned, we take the successive limits $\epsilon_0\searrow 0,\ldots,\epsilon_r\searrow 0$;
in taking this limits, we work modulo terms that are bounded uniformly in the remaining $\epsilon_i$'s. Such terms are {\em any}
terms that contain derivatives of $\phi$. Hence, for the sake of simplicity and bravity, we shall completely drop $\phi$ from the ensuing
coomputations; see also the remark in \cite[p521]{DH} or the proof of \cite[Theorem 4]{BT}.}
\end{remark}

{\bf Proof of Theorem \ref{thm:meli}} We consider the function
\begin{equation}
u(x)=d^{\frac{2mp-k+\epsilon_0}{p}}X_1^{\frac{-1+\epsilon_1}{p}}\ldots
X_r^{\frac{-1+\epsilon_r}{p}},
\label{u}
\end{equation}
where $\epsilon_0,\ldots,\epsilon_r$ are small and positive. 
A standard argument shows that
$u$ lies in the approprite Sobolev space.
We have seen that
\begin{equation}
I_{2m,r-1}[u]=\sum_{0\leq i\leq j\leq r}a_{ij}\Gamma_{ij} + O(1),
\label{111}
\end{equation}
where the coefficients $a_{ij}$ are given by (\ref{aij}) and the $s_i$'s are related to
the $\epsilon_i$'s by (\ref{se}).

We let $\epsilon_0\searrow 0$ in (\ref{rain}). It follows from (\ref{finite}) that all $\Gamma_{ij}$'s
with $i\geq 1$ have finite limits.
As for the remaining terms, applying Lemma \ref{diverge} with $\beta=-3/2$ (for $j=0$)
and with $\beta=-1/2$ (for $j\geq 1$) we obtain respectively
\begin{equation}
\Gamma_{00}\leq c\epsilon_0^{-\frac{5}{2}} \qquad  , \qquad \Gamma_{0j}\leq c\epsilon_0^{-\frac{3}{2}},
\label{xeim}
\end{equation}
where in both cases $c>0$ is independent of all the $\epsilon_i$'s.
Now, we think of the quantities $a_{0j}$ of Lemma \ref{trian} as functions of $\epsilon_0$
and consider $\epsilon_1,\ldots,\epsilon_r$ as small positive parameters.
Using Taylor's theorem we shall expand the coefficient $a_{0j}$ of $\Gamma_{0j}$, $j= 0$ (resp. $j\geq 1$) in powers
of $\epsilon_0$ and (\ref{xeim}) shows that we can discard powers
with exponent $\geq 3$ (resp. $\geq 2$). We shall compute the remaining ones and for this we define
\[
\hat{a}_{0j}(\epsilon_0):=a_{0j}(s_0)=a_{0j}((2mp-k+\epsilon_0)/p)
\]
and denote by $A_{k,0j}$ the coefficient of $\epsilon_0^k$ in $\hat{a}_{0j}(\epsilon_0)$. 
We then have from Lemma \ref{trian}:

{\em -- Constant term in $a_{00}$:} We have 
\begin{eqnarray}
\hat{a}_{00}(\epsilon_0)&=& |\alpha_m(s_0)|^p-|A(2m)| \label{elsa} \\
&=& \left| \prod_{i=0}^{m-1}\Big( \frac{(2m-2i)p-k+\epsilon_0}{p}\Big)
\prod_{j=1}^{m}\Big( \frac{(2m-2j)p+kp-k+\epsilon_0}{p}\Big)\right|^p -|A(2m)|
\nonumber
\end{eqnarray}
and therefore, using (\ref{av1}), $A_{0,00}=\hat{a}_{00}(0)=|\alpha_m(s_0)|^p\Big|_{\epsilon_0=0}-|A(2m)|=0$.

{\em -- Coefficient of $\epsilon_0$ in $a_{00}$:} Differentiating (\ref{elsa}) we obtain
\begin{equation}
\hat{a}_{00}'(\epsilon_0)=\frac{1}{p}a_{00}'(s_0)= |a_m(s_0)|^{p-2}a_m(s_0)a_m'(s_0)
\label{manos}
\end{equation}
and therefore the coefficient is
\[
A_{1,00}=\hat{a}_{00}'(0)=|a_m(s_0)|^{p-2}a_m(s_0)a_m'(s_0)\Big|_{\epsilon_0=0}.
\]

{\em -- Coefficient of $\epsilon_0^2$ in $a_{00}$:} We have from (\ref{manos})
\[
A_{2,00}=\frac{\hat{a}_{00}''(s_0)}{2}\Big|_{\epsilon_0=0}= \frac{1}{2p}
\Big(|a_m(s_0)|^{p-2}a_m(s_0)a_m'(s_0)\Big)'\Big|_{\epsilon_0=0}\; .
\]
Concerning $a_{0j}$, $j\geq 1$, we have $\hat{a}_{0j}(\epsilon_0)=ps_j|\alpha_m(s_0)|^{p-2}\alpha_m(s_0) \alpha_m'(s_0)$
and therefore
\[
\hat{a}_{0j}'(\epsilon_0)=s_j\Big(|\alpha_m(s_0)|^{p-2}\alpha_m(s_0) \alpha_m'(s_0)\Big)'.
\]
Hence:

{\em -- Constant term in $a_{0j}$, $j\geq 1$:} This is
\[
 A_{0,0j}=\hat{a}_{0j}(0)= ps_j|\alpha_m(s_0)|^{p-2}\alpha_m(s_0) \alpha_m'(s_0)\Big|_{\epsilon_0=0}
\]
{\em -- Coefficient of $\epsilon_0$ in $a_{0j}$:} This is
\[
A_{1,0j}=\hat{a}_{0j}'(0)= s_j \Big(|\alpha_m(s_0)|^{p-2}\alpha_m(s_0) \alpha_m'(s_0)\Big)' \Big|_{\epsilon_0=0}.
\]

Now, we observe that $A_{0,0j}=ps_jA_{1,00}=(\epsilon_j-1)A_{1,00}$. Hence $\id$ of Lemma \ref{diverge} implies that
\begin{equation}
A_{1,00}\epsilon_0\Gamma_{00} +\sum_{j=1}^r A_{0,0j}\Gamma_{0j} =O(1)
\label{p1}
\end{equation}
uniformly in $\epsilon_1,\ldots,\epsilon_r$. Similarly, we observe that $A_{1,0j}=2ps_jA_{2,00}=2(-1+\epsilon_j)A_{2,00}$.
Hence, by $\ic$ of Lemma \ref{diverge}, the remaining `bad' terms when combined give
\begin{eqnarray}
&&\hspace{-2cm}A_{2,00}\epsilon_0^2\Gamma_{00} +\epsilon_0\sum_{j=1}^r A_{1,0j}\Gamma_{0j}=\nonumber \\
&=&A_{2,00}\Big( \epsilon_0^2\Gamma_{00} -2\epsilon_0\sum_{j=1}^r (1-\epsilon_j)\Gamma_{0j}\Big)\label{hh}\\
&=&A_{2,00}\Big(
\sum_{i=1}^r(\epsilon_i-\epsilon_i^2)\Gamma_{ii} -\sum_{1\leq i<j\leq r}(1-\epsilon_j)(1-2\epsilon_i)\Gamma_{ij}
\Big)+O(1),
\nonumber
\end{eqnarray}
uniformly in $\epsilon_1,\ldots,\epsilon_r$. Note that the right-hand side of (\ref{hh}) has a finite limit as
$\epsilon_0\searrow 0$.
From (\ref{111}), (\ref{p1}) and (\ref{hh}) we conclude that,
after letting $\epsilon_0\searrow 0$, we are left with
\begin{eqnarray}
&&\qquad I_{2m,r-1}[u] \nonumber \\
&=&\sum_{i=1}^r\bigg( a_{ii}+ A_{2,00}(\epsilon_i-\epsilon_i^2)\bigg)\Gamma_{ii}
+\sum_{1\leq i< j\leq r}\bigg( a_{ij} -A_{2,00}(1-\epsilon_j)(1-2\epsilon_i)  \bigg)\Gamma_{ij}
 +O(1) \nonumber \\
&=:&\sum_{i=1}^rb_{ii}\Gamma_{ii} + \sum_{1\leq i<j\leq r}b_{ij}\Gamma_{ij}+O(1)\; , \quad\quad (\epsilon_0=0),
\label{alv}
\end{eqnarray}
where the $O(1)$ is uniform in $\epsilon_1,\ldots,\epsilon_r$.

We next let $\epsilon_1\searrow 0$ in (\ref{alv}). It follows from (\ref{finite}) that all the $\Gamma_{ij}$'s have
finite limits, except those with $i=1$ which diverge to $+\infty$. For the latter we have
\[
\Gamma_{11}\leq c\epsilon_1^{-\frac{5}{2}} \quad , \qquad
\Gamma_{1j}\leq c\epsilon_1^{-\frac{3}{2}}\;\;, \quad j\geq 2,
\]
by $\ib$ of Lemma \ref{diverge} with $\beta=-3/2$ and $\beta=-1/2$ respectively;
in both cases the constant $c$ is independent of $\epsilon_2,\ldots,\epsilon_r$.
We think of the coefficients $b_{1j}$ as functions -- indeed, polynomials -- of $\epsilon_1$
and we expand these in powers of $\epsilon_1$.
The estimates above on $\Gamma_{1j}$ imply that only the terms $1,\epsilon_1$ and $\epsilon_1^2$
(resp. 1 and $\epsilon_1$) give contributions for $\Gamma_{11}$
(resp. $\Gamma_{1j}$, $j\geq 2$) that do not vanish as $\epsilon_1\searrow 0$.
We shall compute the coefficients of these terms. Our starting point are the relations (cf. (\ref{alv})
\begin{eqnarray}
b_{11}(\epsilon_1)&=&a_{11}(s_0)+A_{2,00}(\epsilon_1-\epsilon_1^2)\nonumber\\
&=&\frac{\epsilon_1-1}{2p}|\alpha_m|^{p-2}\Big( \alpha_m\alpha_m''(p-1+\epsilon_1) 
+(p-1)(\epsilon_1-1)\alpha_m'^2 \Big) \nonumber\\
&&\qquad\qquad +A_{2,00}(\epsilon_1-\epsilon_1^2) -|B(2m)| 
\label{b11}
\end{eqnarray}
and, for $j\geq 2$,
\begin{eqnarray}
b_{1j}(\epsilon_1)&=&a_{1j}(s_0)-A_{2,00}(1-\epsilon_j)(1-2\epsilon_1)\nonumber\\
&=&(\epsilon_j-1)\bigg\{\frac{|\alpha_m|^{p-2}}{2p}\Big[ \alpha_m\alpha_m''(2\epsilon_1+p-2) 
+2(p-1)\alpha_m'^2(\epsilon_1-1)\Big]+\nonumber \\
&&\qquad\qquad\qquad+A_{2,00}(1-2\epsilon_1)\bigg\}.
\label{b1j}
\end{eqnarray}
Hence, denoting by $B_{k,1j}$ the coefficient of $\epsilon_1^k$ in $b_{1j}$, $j\geq 1$, we have:

{\em -- Constant term in $b_{11}$:} This is
\begin{eqnarray*}
B_{0,11}&=&b_{11}(0)\\
&=&a_{11}(s_0)\Big|_{\epsilon_1=0}\\
&=&-\frac{1}{2p}|\alpha_m|^{p-2}\Big( \alpha_m\alpha_m''(p-1) -(p-1)\alpha_m'^2 \Big) -|B(2m)| \\
&=&0,
\end{eqnarray*}
by (\ref{av2}).

{\em -- Coefficient of $\epsilon_1$ in $b_{11}$:} From (\ref{b11}) we obtain
\begin{equation}
b_{11}'(\epsilon_1)=\frac{|\alpha_m|^{p-2}}{2p}\bigg\{ \alpha_m\alpha_m''(2\epsilon_1+p-2)
+(p-1)\alpha_m'^2(2\epsilon_1-2)\bigg\} +A_{2,00}(1-2\epsilon_1)
\label{box}
\end{equation}
and therefore the coefficient is 
\[
B_{1,11}=b_{11}'(0)=\frac{|\alpha_m|^{p-2}}{2p}\Big\{ (p-2)\alpha_m\alpha_m''-2(p-1)\alpha_m'^2 \Big\}+A_{2,00}.
\]
{\em -- Coefficient of $\epsilon_1^2$ in $b_{11}$:} From (\ref{box}),
\[
B_{2,11}=\frac{1}{2}b_{11}''(0)=\frac{|\alpha_m|^{p-2}}{2p}\Big\{\alpha_m\alpha_m'' +(p-1)\alpha_m'^2\Big\}
-A_{2,00}=0.
\]
{\em -- Constant term in $b_{1j}$, $j\geq 2$:} This is
\[
B_{0,1j}=b_{1j}(0)=(\epsilon_j-1)\frac{|\alpha_m|^{p-2}}{2p}\Big[(p-2)\alpha_m\alpha_m''-2(p-1)\alpha_m'^2\Big]
+A_{2,00}.
\]
{\em -- Coefficient of $\epsilon_1$ in $b_{1j}$, $j\geq 2$:} From (\ref{b1j}),
\[
B_{1,1j}=b_{1j}'(0)=(\epsilon_j-1)\Big\{  \frac{|\alpha_m|^{p-2}}{p}(\alpha_m\alpha_m'' +(p-1)\alpha_m'^2)
-2A_{2,00}\Big\}=0.
\]
We obsrerve that $B_{0,1j}=(\epsilon_j-1)B_{1,11}$, $j\geq 2$. Hence $\id$ of Lemma \ref{diverge} gives
\begin{equation}
\epsilon_1 B_{1,11}\Gamma_{11}+\sum_{j=2}^rB_{0,1j}\Gamma_{1j} =O(1),
\label{p3}
\end{equation}
uniformly in $\epsilon_2,\ldots,\epsilon_r$.
Combining (\ref{alv}) and  (\ref{p3}) we conclude that after letting $\epsilon_1\searrow 0$ we are left with
\begin{equation}
I_{2m,r-1}[u]=\sum_{2\leq i\leq j\leq r}b_{ij}\Gamma_{ij} +O(1)\; , \qquad\quad (\epsilon_0=\epsilon_1=0),
\label{alv1}
\end{equation}
uniformly in $\epsilon_2,\ldots,\epsilon_r$. Note that we have the
same coefficients $b_{ij}$ as in (\ref{alv}), unlike the case where the limit $\epsilon_0\searrow 0$ was taken,
where we passed from the coefficients $a_{ij}$ to the coefficients $b_{ij}$.

We proceed in this way. At the $i$th step we denote by $B_{k,ij}$ the coefficient of $\epsilon_i^k$ in $b_{ij}$, $j\geq i$, and observe
that
\[
B_{0,ij}=(\epsilon_j-1)B_{1,ii}\;\quad  , \qquad\quad B_{2,ii}=B_{1,ij}=0\; , \qquad j\geq i+1.
\]
Hence $\id$ of Lemma \ref{diverge} implies the cancelation (modulo uniformly bounded terms)
of all terms that, separetely, diverge as $\epsilon_i\searrow 0$. Eventually, after letting $\epsilon_{r-1}\searrow 0$,
we arrive at
\begin{equation}
I_{2m,r-1}[u]=b_{rr}\Gamma_{rr}+O(1) \; , \quad (\epsilon_0=\epsilon_1=\ldots =\epsilon_{r-1}=0).
\label{alvr}
\end{equation}
Since
\[
\int_{\Omega}\frac{|u|^p}{d^{2mp}}X_1^2\ldots X_r^2dx=\Gamma_{rr}
\]
and $\lim_{\epsilon_r\searrow 0}\Gamma_{rr}= +\infty$ (cf (\ref{diverge})) we conclude that
\begin{eqnarray*}
\inf_{C^{\infty}_c(\Omega\setminus K)}\frac{I_{2m,r-1}[v]}{\int_{\Omega}\frac{|v|^p}{d^{2mp}}X_1^2\ldots X_r^2dx}
&\leq&\lim_{\epsilon_r\searrow 0}\frac{b_{rr}\Gamma_{rr}+O(1)}{\Gamma_{rr}}\\
&=& \lim_{\epsilon_r\searrow 0}b_{rr}\\
&=& \lim_{\epsilon_r \searrow 0}a_{rr}\\
&=&\lim_{\epsilon_r\searrow 0} \frac{ps_r}{2}|\alpha_m|^{p-2}\Big( \alpha_m\alpha_m''(s_r+1) +(p-1)\alpha_m'^2s_r\Big) \\
&=&\frac{p-1}{2p}|\alpha_m|^{p-2}(\alpha_m'^2-\alpha_m\alpha_m'')\\
\mbox{(by (\ref{av2}))}\qquad &=&|B(2m)|.
\end{eqnarray*}
This proves part $\ib$ of the theorem. Part $\ia$ follows by slightly modifying the above argument; we omit the details.
$\hfill //$

\



\end{document}